\documentclass{amsart}
\usepackage{amsthm, amsmath}
\usepackage{epsf}
\usepackage{colordvi}
\usepackage{latexsym}
\usepackage{amssymb}
\newcommand{\cellsize}{13}
\newlength{\cellsz} 
\newsavebox{\cell}
\sbox{\cell}{\begin{picture}(\cellsize,\cellsize)
\put(0,0){\line(1,0){\cellsize}}
\put(0,0){\line(0,1){\cellsize}}
\put(\cellsize,0){\line(0,1){\cellsize}}
\put(0,\cellsize){\line(1,0){\cellsize}}
\end{picture}}
\newcommand\cellify[1]{\def\thearg{#1}\def\nothing{}%
\ifx\thearg\nothing
\vrule width0pt height\cellsz depth0pt\else
\hbox to 0pt{\usebox{\cell} \hss}\fi%
\vbox to \cellsz{
\vss
\hbox to \cellsz{\hss$#1$\hss}
\vss}}
\newcommand\tableau[1]{\vtop{\let\\\cr
\baselineskip -16000pt \lineskiplimit 16000pt \lineskip 0pt
\ialign{&\cellify{##}\cr#1\crcr}}}
\theoremstyle{plain}
\newtheorem{thm}{Theorem}

\newtheorem{remark}[thm]{Remark}
\newenvironment{rem}{\begin{remark}\rm}{\end{remark}}
\newtheorem{example}[thm]{Example}

\newenvironment{ex}{\begin{example}\rm}{\end{example}}

\newcommand{\Gr}{\operatorname{Gr}}
\newcommand{\Fl}{\operatorname{Fl}}

\newcommand{\groth}{{\mathfrak G}}

\newcommand{\calO}{{\mathcal O}}

\newcommand{\calS}{{\mathcal S}}

\newcommand{\C}{{\mathbb C}}

\newcommand{\ignore}[1]{}
\begin{document}
\bibliographystyle{amsplain}

\title[Quiver coefficients are Schubert structure Constants]
{Quiver coefficients are Schubert \\ structure constants}  

\author{Anders Skovsted Buch}
\address{Matematisk Institut\\
         Aarhus Universitet\\
         Ny Munkegade, 8000 {\AA}rhus C, Denmark}
\email{abuch@imf.au.dk}

\author{Frank Sottile}
\address{Department of Mathematics\\
        University of Massachusetts\\
        Amherst, MA \ 01003\\
        USA}
\email{sottile@math.umass.edu}

\author{Alexander Yong}
\address{Department of Mathematics\\
         University of California\\
         Berkeley, CA, 94720, USA}
\email{ayong@math.berkeley.edu}

\date{\today}
\thanks{Second author supported in part by NSF CAREER grant DMS-0134860}
\maketitle

\section{Introduction}

Buch and Fulton established a formula for the cohomology class of a
quiver variety~\cite{BF99}, which Buch later extended to
$K$-theory~\cite{Bu01}. The $K$-theory formula expresses a quiver
class as an integral linear combination of products of stable
Grothendieck polynomials.  The {\em quiver coefficients\/} of this
linear combination were conjectured to be nonnegative in cohomology
and to alternate in sign in $K$-theory.  These conjectures were
recently proved by Knutson, Miller, and Shimozono~\cite{KMS},
Buch~\cite{buch:altsign}, and Miller~\cite{miller:altsign}.

Buch, Kresch, Tamvakis, and Yong \cite{BKTY:schquiv, BKTY:groth} gave
combinatorial formulas for the {\em decomposition coefficients} expressing a
Grothendieck polynomial as an integer linear combination of products of
stable Grothendieck polynomials.  In particular, it was proved that
the decomposition coefficients alternate in sign.

Alternation in sign also occurs in 
the Schubert calculus of the flag variety; Brion~\cite{brion} proved
that the $K$-theory Schubert structure constants alternate in sign.

We give natural and explicit equalities between the three aforementioned
integers.
Our argument uses results of Bergeron and Sottile~\cite{BS98} and Lenart,
Robinson, and Sottile~\cite{LRS} who earlier established a connection between the
decomposition coefficients and the Schubert structure constants.
The other main ingredient is the ratio
formula of Knutson, Miller, and Shimozono~\cite{KMS}, as well as some
combinatorial properties of this formula~\cite{buch:altsign,KMS}.
We also give a direct argument that the decomposition coefficients have
alternating signs, based on Brion's theorem, which  
then implies that quiver coefficients have alternating signs.
A consequence of our theorem is that formulas for the other numbers give
formulas for the quiver coefficients; we give examples in the last
section.

\section{The Main Result}

Quiver coefficients $c_{\underline \mu}(r)$ are defined for a set of
rank conditions $r=\{r_{i,j}\}$ for $0\leq~i\leq~j\leq n$
and a sequence of partitions ${\underline\mu} = (\mu_1,\mu_2, \dotsc,
\mu_n)$, where $\mu_i$ fits in a $r_{i-1}\times r_i$ rectangle (for
convenience, set $r_i:=r_{i,i}$).  The Grothendieck polynomial
$\groth_w$ for a permutation $w$ represents the class of the structure
sheaf of the corresponding Schubert variety in the Grothendieck ring
of the flag variety~\cite{LS:groth}.  These form a ${\mathbb
  Z}$-linear basis for the Grothendieck ring. The {\em Schubert
  structure constants} are the unique integers $C_{u,v}^{w}$ such that
\[ \groth_u\cdot\groth_v\ =\ \sum_w C^w_{u,v}\ \groth_w \,. \]
Brion's theorem \cite[Thm.~1]{brion} proves that
$(-1)^{\ell(wuv)}C^w_{u,v}\geq 0$.

For a partition $\lambda$, let $w(\lambda,k)$ denote the Grassmannian
permutation for $\lambda$ with descent at $k$.  It is given by
$w(\lambda,k)(i)=i+\lambda_{k+1-i}$ for $1\leq i\leq k$ and
$w(\lambda,k)(i)<w(\lambda,k)(i+1)$ for $i\neq k$.  The Grothendieck
polynomial $\groth_{w(\lambda,k)}(y_1,y_2,\dotsc)$ is symmetric
in the variables $y_1,\dotsc,y_k$ and is independent of
$y_i$ for $i>k$.  Thus we can write $G_\lambda(y_1,\dotsc,y_k)$ for
this symmetric Grothendieck polynomial, without ambiguity.

Suppose that $v$ is a permutation whose descents occur at positions
in $\{r_0, r_0+r_1, \dotsc,r_0+\dotsb+r_{n-1}\}$.  If
$x^i = (x^i_1,x^i_2,\dotsc,x^i_{r_i})$ is a set of $r_i$ variables,
then the Grothendieck polynomial $\groth_v(x^0,x^1,\dotsc,x^n)$ is
separately symmetric in each set of variables $x^i$ and is independent
of $x^n$.  As the $G_\lambda$ form a basis for all symmetric
polynomials, there are integer {\it decomposition coefficients}
$b_{\underline \mu}(v)$ defined by
\begin{equation}\label{E:decomposition}
   \groth_v(x^0,x^1,\dotsc,x^n)\ =\ \sum_{\underline\mu}
   b_{\underline \mu}(v)\, G_{\mu_1}(x^0)\, G_{\mu_2}(x^1)\dotsb 
   G_{\mu_n}(x^{n-1})\,.
\end{equation}
Formulas for these coefficients given in~\cite{BKTY:groth} show that
$(-1)^{\sum |\mu_i| - \ell(v)} b_{\underline \mu}(v) \geq 0$.  In
Remark~\ref{R:simple_geometry} below, we give a simple geometric
argument that accounts for this via Brion's theorem.

In order to state our main result, we will need some notation and
terminology.  For a set of rank conditions $r$, let $d'_i :=
r_i+\dotsb+r_{n-1}$ and $d_i:=d'_i+r_n$.  Also let $R_i =
(d_{i+1})^{r_{i-1}}$ be the rectangular partition with $r_{i-1}$ rows
and $d_{i+1}$ columns.  For a sequence of partitions ${\underline \mu}
= (\mu_1,\mu_2,\dotsc,\mu_n)$ let $\widetilde{\mu}_i$ be the result of
attaching $\mu_i$ to the right side of $R_i$, and let
${\underline{\widetilde{\mu}}} = ({\widetilde{\mu}_1}, \ldots,{\widetilde{\mu}_n})$
denote the sequence of these partitions.
\begin{equation}\label{mutilde}
  \raisebox{-24pt}{
  \begin{picture}(200,40)(-40,0)
  \put(-40,24){$\widetilde{\mu}_i\ :=$}
  \put(22,10){\epsffile{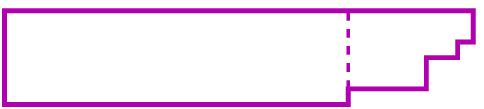}}
  \put(0,24){$r_{i-1}$}
  \put(60,0){$d_{i+1}$}
  \put(130,25){$\mu_i$}
  \put(70,25){$R_i$}
  \end{picture}}
\end{equation}
Let $\rho$ be a concatenation of rectangles of sizes $d'_i\times d_i$
for $i=1,2,\dotsc,n{-}1$.  This is the shaded part of the partition
shown in~\eqref{E:rho}.  Finally, we let $\rho({\underline\mu})$ be
the partition obtained when $\widetilde{\mu}_i$ is placed under the $i$th
rectangle of $\rho$ and the result is concatenated with $\mu_n$.
This has $d:=d'_0$ rows and is shown below.
\begin{equation}\label{E:rho}
\raisebox{-45pt}{
  \setlength{\unitlength}{0.8pt}
  \begin{picture}(360,125)(-17,0)
  \put(0,0){\epsfbox{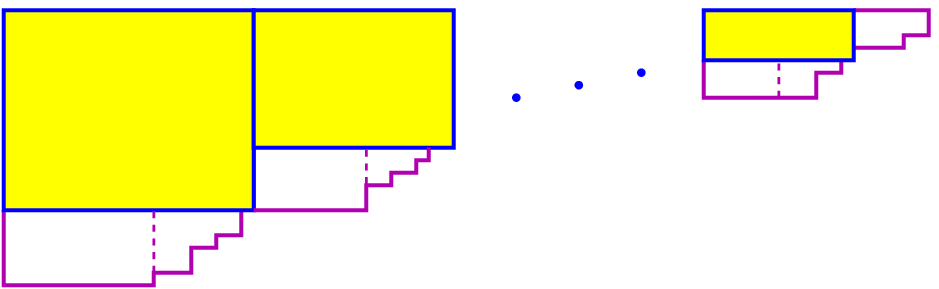}} 
  \put(-17, 65){$d'_1$}
  \put( 35,110){$d_1$}
  \put( 80,  3){$\widetilde{\mu}_1$}

  \put(168, 78){$d'_2$}
  \put(120,110){$d_2$}
  \put(143, 25){$\widetilde{\mu}_2$}

  \put(218, 90){$d'_{n-1}$}
  \put(270,110){$d_{n-1}$}
  \put(270, 56){$\widetilde{\mu}_{n-1}$}

  \put(326,79){$\mu_n$}

  \end{picture}
}
\end{equation}

The Zelevinsky permutation $v(r)$ encodes the rank conditions $r$; we
refer to \cite[Def.~1.7]{KMS} for the precise definition.  We note
that the descents of $v(r)$ occur at positions in 
$\{r_0, r_0+r_1, \dotsc,r_0+\dotsb+r_{n-1}\}$.

\begin{thm}\label{thm:main}
  Let $r$ be a set of rank conditions and ${\underline \mu}=
  (\mu_1,\ldots,\mu_n)$ be a sequence of partitions with $\mu_i$ a
  subset of the $r_{i-1}\times r_{i}$ rectangle.  Then the following
  numbers are equal:

\begin{itemize}
 \item[(I)] the quiver coefficient $c_{\underline\mu}(r)$;
 \item[(II)] the decomposition coefficient $b_{\widetilde{\underline\mu}}(v(r))$;
 \item[(III)] the Schubert structure constant 
               $C^{w(\rho({\underline\mu}),d)}_{v(r),\,w(\rho,d)}$.
\end{itemize}
\end{thm}

We remark that we follow the notation for quiver coefficients used in
\cite{KMS}, which conjugates all partitions compared to the notation
used in \cite{BF99,Bu01,buch:altsign}.

\begin{proof}
The ratio formula of Knutson, Miller, and Shimozono expresses a quiver
class as a quotient of two Grothendieck polynomials \cite[Thm.~2.7]{KMS}.
We need the following identity which is derived from the ratio formula
in \cite[\S 7]{buch:altsign}:
\begin{equation} \label{E:I.II}
  \frac{\groth_{v(r)}(x^0,x^1,\dotsc,x^n)}
  {\groth_{v(\varepsilon)}(x^0,x^1,\dotsc,x^n)}\ =\ 
  \sum_{\underline\mu} c_{\underline\mu}(r)\, G_{\mu_1}(x^0) \,
  G_{\mu_2}(x^1)\dotsb G_{\mu_n}(x^{n-1}) \,.
\end{equation}
Here $\epsilon$ denotes the maximal rank conditions given by
$\epsilon_{ij} = \min\{r_i,\dots,r_j\}$.  Notice that the cohomology
version of (\ref{E:I.II}) follows from Thm.~7.10 and Prop.~7.13 of
\cite{KMS}.  The denominator of (\ref{E:I.II}) is the monomial
\[ \groth_{v(\epsilon)}(x^0,x^1,\dots,x^n) \ =\  
   \prod_{i=1}^{n-1} G_{R_i}(x^{i-1}) \ =\ 
   \prod_{i=1}^{n-1} (x^{i-1}_1 x^{i-1}_2 
     \cdots x^{i-1}_{r_{i-1}})^{d_{i+1}} \,.
\]
By using that $G_{R_i}(x^{i-1})\, G_{\mu_i}(x^{i-1}) = 
   G_{\widetilde{\mu}_i}(x^{i-1})$ \cite[Cor.~6.5]{buch:LR} we deduce that
\[
  \groth_{v(r)}(x^0,x^1,\dotsc,x^n) =
  \sum_{\underline\mu} c_{\underline\mu}(r)\, G_{\widetilde{\mu}_1}(x^0) \,
  G_{\widetilde{\mu}_2}(x^1)\dotsb G_{\widetilde{\mu}_n}(x^{n-1}) \,.
\]
This proves the equivalence of (I) and (II).  Now since the last
descent of $v(r)$ occurs before position $d$, it follows from
\cite[Thm.~9.7]{LRS} that
\begin{equation}\label{E:9.7}
   \groth_{v(r)}(x^0,x^1,\dotsc,x^n)\ =\ 
   \sum_{u_i\in\calS_{d_{i-1}}} \ 
   C^{(u_1\times\dotsb\times u_n)\cdot w(\rho,d)}_{v,\,w(\rho,d)}
   \groth_{u_1}(x^0)\dotsb\groth_{u_n}(x^{n-1}) \,,
\end{equation}
where $\calS_{d_{i-1}}$ is the symmetric group on $d_{i-1}$ elements,
and $u_1 \times \cdots \times u_n \in \calS_{d_0+\dots+d_{i-1}}$ is
the Cartesian product of the permutations $u_i$.  Because the left
hand side is symmetric in each set of variables $x^i$, it follows that
all permutations $u_i$ that occur with non-zero coefficient must be
Grassmannian with descent at position $r_{i-1}$, so $u_i =
w(\lambda_i,r_{i-1})$ for a partition $\lambda_i$.  Furthermore, since
$G_{R_i}(x^{i-1})$ divides the left hand side of (\ref{E:9.7}), each
partition $\lambda_i$ must have the form $\widetilde{\mu}_i$ for a
partition $\mu_i$.  Using that
\[ w(\rho({\underline\mu}),d) \ =\ 
   \bigl( w(\widetilde{\mu}_1,r_0)\times \dotsb \times  
          w(\widetilde{\mu}_n,r_{n-1})\bigr)\cdot w(\rho,d) \,,
\]
we deduce that
\begin{equation}\label{E:subs}
  \groth_{v(r)}(x^0,x^1,\dotsc,x^n)\ =\ 
  \sum_{\underline\mu} C^{w(\rho({\underline\mu}),d)}_{v(r), w(\rho,d)} \;
  G_{\widetilde{\mu}_1}(x^0)\, G_{\widetilde{\mu}_2}(x^1) \dotsb
  G_{\widetilde{\mu}_n}(x^{n-1})\,.
\end{equation}
This proves the equivalence of (II) with (III).
\end{proof}

\begin{ex}\label{EX:one}
Suppose that $n=3$ and let $r$ be the set of rank conditions
\[ 
  \begin{matrix}
    r_{00} && r_{11} && r_{22} && r_{33} \\
    & r_{01} && r_{12} && r_{23} \\
    && r_{02} && r_{13} \\
    &&& r_{03}
  \end{matrix}
   \qquad=\qquad
  \begin{matrix}
    1 & & 4 & & 3 & & 3 \\
     & 1 & & 2 & & 2 \\
      & & 1 & & 1 \\
       & & & 0
  \end{matrix}
\]
Here, $(d'_0,d'_1,d'_2)=(8,7,3)$, $(d_0,d_1,d_2,d_3)=(11,10,6,3)$, and
the Zelevinsky permutation $v(r)$ is
$(7,\,4,5,8,9,\,1,2,11,\,3,6,10)$.  The partition $\rho$ is the
concatenation of a $7\times 10$ rectangle with a $3\times 6$
rectangle, and so equals $(16,16,16, 10,10,10,10)$.  This is the
shaded part of the partition shown in~\eqref{example}.

Let ${\underline\mu}=(\emptyset, \ (2,1,1), \ (1))$ be a sequence of
partitions.  Then $\underline{\widetilde{\mu}} = ((6), \ (5,4,4,3), \linebreak (1))$
and the partition $\rho(\underline\mu)$ is $(17,16,16, 15,14,14,13,
6)$, which is illustrated below.
\begin{equation}\label{example}
  \raisebox{-42pt}{\epsfbox{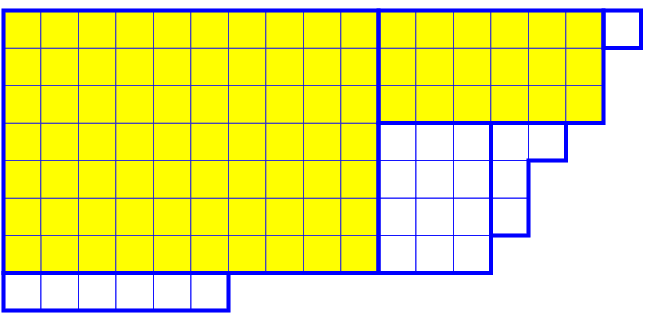}} 
\end{equation}

In Example 1 of~\cite{Bu01}, the quiver coefficient
$c_{\underline\mu}(r)$ was computed to be 1.  Due to our different
conventions, the corresponding term there is $1 \otimes
s_{\epsffile{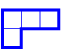}}\otimes s_{\epsffile{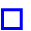}}$,
which is indexed by the sequence of partitions $(\emptyset, (3,1), (1))$.
Later, we will use formulas from~\cite{BKTY:schquiv}
and~\cite{Kogan:schur_schub} to calculate that the corresponding
decomposition coefficients and Schubert structure constants are both~1.
\end{ex}

\begin{rem} \label{R:simple_geometry}
  The alternating signs of the decomposition coefficients can be explained
  geometrically using Brion's theorem as follows.  Let $v$ be a
  permutation whose descents occur at positions
  in $\{r_0, r_0+r_1, \dotsc,r_0+\dotsb+r_{n-1}\}$. 
  Fix a large integer $N$ and let $Y =
  \Fl(r_0,r_0+r_1,\dots,r_0+\cdots+r_{n-1} ; \C^{nN})$ be the variety
  of partial flags of the indicated type in $\C^{nN}$.  
  Then $\groth_v$ represents the class $[\calO_{Y_v}]$ of the structure
  sheaf of the Schubert variety $Y_v$ of $Y$ \cite{LS:groth}.  The product of
  Grassmannians
\[ 
    X = \Gr(r_0,N) \times \Gr(r_1,N) \times \dots \times \Gr(r_{n-1},N)
\]
can be identified with a subvariety of $Y$ by mapping each point
$(V_1,\dots,V_n) \in X$ to the partial flag $V_1 \subset V_1\oplus V_2
\subset \cdots \subset V_1 \oplus \cdots \oplus V_n$ in $Y$.  If we
let $x^i = \{x^i_1,\dots,x^i_{r_i}\}$ denote the $K$-theoretic Chern
roots of the dual of the tautological subbundle corresponding to the
$i$th factor of $X$, then the specialization $\groth(x^0,x^1,...,x^n)$
is obtained by restricting the class $[\calO_{Y_v}]$ to the
Grothendieck ring of $X$.  When the Schubert variety $Y_v$ is in
general position, it follows from \cite[Lemma~2]{brion} that the
subvariety $X_v := X \cap Y_v$ of $X$ has rational singularities.  We
furthermore obtain that
\[ [\calO_{X_v}] = \groth_v(x^0,x^1,\dotsc,x^n) = 
   \sum b_{\underline \mu}(v)\, G_{\mu_1}(x^0)\,
   G_{\mu_2}(x^1)\dotsb G_{\mu_n}(x^{n-1})
\]
in the Grothendieck ring of $X$.
Now since the right hand side of this identity is a linear combination
of Schubert structure sheaves on the generalized flag variety $X$, we
conclude from \cite[Thm.~1]{brion} that $(-1)^{\sum |\mu_i| - \ell(v)}
b_{\underline \mu}(v) \geq 0$, as required.

In~\cite{BS98} similar geometry is used to study the
restriction of Schubert classes $[Y_v]$ in a full flag manifold $Y$ to
products $X$ of flag manifolds embedded in $Y$ in a similar fashion to that
given here.
There, the subvariety $X_v$ is identified as an intersection of Schubert
varieties in $Y$.
This is used to identify the decomposition coefficients as particular
Schubert structure constants, for cohomology.
The analogous identification of $K$-theoretic decomposition coefficients with
$K$-theoretic Schubert structure constants is accomplished in~\cite{LRS}.
\end{rem}

\section{Alternative formulas for quiver coefficients}

By Theorem~\ref{thm:main}, 
formulas for the decomposition coefficients and Schubert structure
constants give alternative formulas for the quiver coefficients. 
We give three examples of this for the
cohomology quiver coefficients (which are indexed by sequences of
partitions such that $\sum |\mu_i|$ equals the expected codimension of the
quiver variety, see~\cite{BF99}).

The formulas of~\cite{BKTY:schquiv,BKTY:groth} for the decomposition
coefficients give formulas for quiver
coefficients in cohomology and in $K$-theory.  
We state this formula for the cohomology decomposition coefficients.  
Suppose that $v$ is a
permutation whose descents occur at positions
 in $\{r_0, r_0+r_1, \dotsc,r_0+\dotsb+r_{n-1}\}$.
 The decomposition coefficient
$b_{\underline \mu}(v)$ is equal to the number of sequences of
semistandard tableaux $(T_1,\ldots,T_n)$ such that
\begin{itemize}
\item[(i)]   The shape of $T_i$ is the conjugate (matrix transpose) of the
  partition $\mu_i$,
\item[(ii)] The entries of $T_i$ are strictly larger than
  $r_0+\dotsb+r_{i-2}$, and
\item[(iii)] Concatenating the bottom-up, left-to-right column reading 
words of the tableaux $T_1, T_2, \dotsc, T_n$ gives a reduced word for $v$.
\end{itemize}

The quiver coefficient computed in Example~\ref{EX:one} corresponds to the 
sequence of tableaux
\[ \tableau{1\\ 2\\ 3\\ 4\\ 5\\ 6} \ \ \ \ \ 
   \tableau{2& 3& 4& 5\\ 3& 4& 5& 6\\ 4& 5& 6& 7\\ 7& 8& 9\\ 10} \ \ \ \ \ 
   \tableau{8}
\]
which encodes the following reduced word:
\[v(r)=s_{6}s_{5}s_{4}s_{3}s_{2}s_{1}s_{10}s_{7}s_{4}s_{3}s_{2}s_{8}s_{5}s_{4}
s_{3}s_{9}s_{6}s_{5}s_{4}s_{7}s_{6}s_{5}s_{8} \,. \]

As another example, 
Knutson~\cite{Knutson:action} gives an algorithm to compute any
cohomology Schubert structure constant $C^w_{u,v}$, i.e., when 
$\ell(u)+\ell(v)=\ell(w)$, where $\ell(u)$ is the length of a permutation
$u$. While this involves signs in general, there
are no signs when $w$ is a Grassmannian permutation.
Thus, this gives a subtraction-free algorithm to compute the 
quiver coefficient 
$c_{\underline\mu}(r)=C^{w(\rho({\underline\mu}),d)}_{v(r),w(\rho,d)}$.

Finally, Kogan~\cite{Kogan:schur_schub} gives a generalization of the
Littlewood-Richardson rule for the cohomology Schubert structure
constants $C^u_{v,w(\lambda,d)}$, when $v$ has no descents after position $d$.
We give a mild reformulation of his formula in terms of 
chains in the Bruhat order from $v$ to $u$ that give
valid reading words for tableaux of shape $\lambda$.

A {\em saturated chain} $\gamma$ in the $d$-Bruhat order is a sequence
of permutations
\[
  \gamma\ \colon\ 
   v\ =\ v_0\ \longrightarrow\ v_1\ \longrightarrow\ v_2\  
    \longrightarrow\ \dotsb\ \longrightarrow\ v_t\ =\ u\,,
\]
where $\ell(v_i)=\ell(v)+i$ and $v_{i-1}^{-1}v_{i}$ is a transposition
$(j_i,k_i)$ with $j_i\leq d<k_i$ for each $i=1,\dotsc,t$.
The {\it word} of such a chain $\gamma$ is the sequence of integers
\[
  v_1(k_1),\, v_2(k_2),\,\dotsc\,,\, v_t(k_t)\,.
\]
Kogan's formula~\cite[Theorem 2.4]{Kogan:schur_schub} asserts that if
$v$ has no descents after position $d$, then $C^u_{v,w(\lambda,d)}$ is
equal to the number of saturated chains in the $d$-Bruhat order from
$v$ to $u$ whose word is the left-to-right, bottom-up row reading word
of a tableau of shape $\lambda$.

For the quiver coefficient of Example~\ref{EX:one}, we have $d=8$ and
there is exactly one saturated chain in the $d$-Bruhat order that goes
from $v(r)$ to $w(\rho({\underline\mu}))$ and whose word is a reading
word for a semistandard tableau of shape $\rho$.

To describe it, we define a chain in the $d$-Bruhat
order to be {\em increasing} if its word is an increasing sequence.
If there is an increasing chain from $v$ to $u$, then it is unique, the
permutation $vu^{-1}$ is the product of disjoint cycles where the numbers
decrease in each cycle, and the partition of the numbers by the cycles they lie
in is non-crossing~\cite[p.~655]{BS02}.
The desired chain has length 88 and it is the concatenation
of increasing chains of lengths $10, 10, 10, 10, 16, 16, 16$, respectively.
Below, we display each increasing chain on a separate line.
Each product of cycles on a given line is $v_j v_i^{-1}$, 
where the increasing subchain for that line is from $v_i$ to $v_j$.
\[
  \begin{array}{l}
   (21, 20, 19, 18, 17, 16, 15, 14, 13, 12, 11)\\
   (19, 18, 17, 16, 15, 14, 13, 12, 11, 10,  9)\\
   (18, 17, 16, 15, 14, 13, 12, 11, 10, 9,  8)\\
   (16, 15, 14, 13, 12, 11, 10, 9,  8, 6, 5)\\
   (25, 24, 23, 22, 21)(20, 19)(17, 16)(15,14,13,12,11,10, 9, 8, 6, 5, 4)\\
   (23, 22, 21, 19, 16, 14, 13, 12, 11, 10, 9, 8, 6, 5, 4, 3, 3,  2)\\
   (22, 21, 19, 16, 14, 13, 12, 11, 10, 9, 8, 6, 5, 4, 3, 2, 1)
  \end{array}
\]
This chain corresponds to the following semistandard tableau of shape $\rho$.

\[ \tableau{
    1&  2&  3&  4&  5&  6&  8&  9& 10& 11& 12& 13& 14& 16& 19& 21\\
    2&  3&  4&  5&  6&  8&  9& 10& 11& 12& 13& 14& 16& 19& 21& 22\\
    4&  5&  6&  8&  9& 10& 11& 12& 13& 14& 16& 19& 21& 22& 23& 24\\
    5&  6&  8&  9& 10& 11& 12& 13& 14& 15\\
    8&  9& 10& 11& 12& 13& 14& 15& 16& 17\\
    9& 10& 11& 12& 13& 14& 15& 16& 17& 18\\
   11& 12& 13& 14& 15& 16& 17& 18& 19& 20
}\]

%
%
%

\begin{rem}
  It would be interesting to give bijections between the 
  formulas discussed here
  for the quiver coefficients (in cohomology) with 
  those given in~\cite{KMS}.
\end{rem}



\begin{thebibliography}{99}
\providecommand{\bysame}{\leavevmode\hbox to3em{\hrulefill}\thinspace}
\providecommand{\MR}{\relax\ifhmode\unskip\space\fi MR }
\providecommand{\MRhref}[2]{%
  \href{http://www.ams.org/mathscinet-getitem?mr=#1}{#2}
}
\providecommand{\href}[2]{#2}

\bibitem{BS98}
N.~Bergeron and F.~Sottile, \emph{Schubert polynomials, the {B}ruhat
  order, and the geometry of flag manifolds}, Duke Math. J. \textbf{95} (1998),
  no.~2, 373--423. \MR{2000d:05127}

\bibitem{BS02}
\bysame, \emph{Skew {S}chubert functions and the {P}ieri formula for flag
  manifolds}, Trans. Amer. Math. Soc. \textbf{354} (2002), no.~2, 651--673
  (electronic). \MR{1 862 562}

\bibitem{billey.braden:patterns}
S.~Billey and T.~Braden, \emph{Lower bounds for Kazhdan-Luztig polynomials from
patterns}, to appear in Transformation Groups, 2003.

\bibitem{brion}
M.~Brion, \emph{Positivity in the Grothendieck group of complex flag
varieties}, J.~Algebra {\bf 258} (2002), no.~1, 137--159.

\bibitem{buch:altsign}
A.~S.~Buch, \emph{Alternating signs of quiver coefficients},
preprint, 2003.

\bibitem{Bu01}
\bysame, \emph{Grothendieck classes of quiver varieties},
Duke Math. J. {\bf 115} (2002), no.~1, 75--103.

\bibitem{buch:LR}
\bysame, \emph{A Littlewood-Richardson rule for the $K$-theory of
Grassmannians}, Acta Math.~{\bf 189} (2002), 37--78.

\bibitem{BF99}
A.~S.~Buch and W.~Fulton, \emph{Chern class formulas
for quiver varieties}, Invent. Math. {\bf 135} (1999), no.~3, 665--687.

\bibitem{BKTY:schquiv}
A.~S.~Buch, A.~Kresch, H.~Tamvakis and A.~Yong,
\emph{Schubert polynomials and quiver formulas}, to appear in
Duke Math.~J., 2002.

\bibitem{BKTY:groth}
\bysame,
\emph{Grothendieck polynomials and quiver formulas}, preprint, 2003.

\bibitem{fomin.kirillov:groth}
S.~Fomin and A.~N.~Kirillov, {\em Grothendieck polynomials and the Yang-Baxter
equation}, Proceedings of the 6th Inter.~Conf.~on Formal Power Series and Algebraic Combinatorics, DIMACS (1994), 183--190.

\bibitem{fomin.kirillov:yang}
\bysame, {\em The Yang-Baxter equation, symmetric functions, and Schubert polynomials}, Proceedings of the 5th Conference on
Formal Power Series and Algebraic Combinatorics (Florence, 1993), Discrete Math. {\bf 153} (1996), no.~1-3,
123--143.

\bibitem{Knutson:action}
A.~Knutson, \emph{A Schubert calculus recurrence from the noncomplex
$W$-action on $G/B$}, preprint, 2003.

\bibitem{KMS}
A.~Knutson, E.~Miller and M.~Shimozono, \emph{Four positive
formulae for type $A$ quiver polynomials}, preprint, 2003.

\bibitem{Kogan:schur_schub}
M.~Kogan, \emph{RC-graphs and a generalized Littlewood-Richardson rule}, 
Internat. Math. Res. Notices 2001, no. 15, 765--782. 

\bibitem{LS82a}
Alain Lascoux and M.-P.~Sch{\"u}tzenberger, \emph{Polyn\^omes de
  {S}chubert}, C. R. Acad. Sci. Paris S\'er. I Math. \textbf{294} (1982),
  no.~13, 447--450. \MR{83e:14039}

\bibitem{LS:groth}
\bysame, {\em Structure de Hopf de
l'anneau de cohomologie et de l'anneau de Grothendieck d'une
vari\'{e}t\'{e} de drapeaux}, C.R.~Acad.~Sci.~Paris S\'{e}r.~I Math.
{\bf 295} (1982), no.~11, 629--633.

\bibitem{LRS}
C.~Lenart, S.~Robinson and F.~Sottile, {\em Grothendieck polynomials via
permutation patterns and chains in the Bruhat order}, 2003.

\bibitem{miller:altsign}
E.~Miller, {\em Alternating formulae for $K$-theoretic quiver
polynomials}, 2003.

\end{thebibliography}
\end{document}